\documentclass[12pt,leqno]{amsart}  
\usepackage{amsmath,amstext,amsthm,amssymb}
\usepackage[colorlinks,citecolor=red,pagebackref,hypertexnames=false]{hyperref}
\usepackage[backrefs]{amsrefs}
\allowdisplaybreaks
\swapnumbers
\theoremstyle{plain}

\newtheorem{lemma}[equation]{Lemma} 
\newtheorem{proposition}[equation]{Proposition} 
\newtheorem{theorem}[equation]{Theorem} 
\newtheorem{corollary}[equation]{Corollary}

\def\norm#1.#2.{\lVert#1\rVert_{#2}}
\def\Norm#1.#2.{\bigl\lVert#1\bigr\rVert_{#2}}
\def\NOrm#1.#2.{\Bigl\lVert#1\Bigr\rVert_{#2}}
\def\NORm#1.#2.{\biggl\lVert#1\biggr\rVert_{#2}}
\def\NORM#1.#2.{\Biggl\lVert#1\Biggr\rVert_{#2}}

\def\eqdef{\stackrel{\mathrm{def}}{{}={}}}

\def\ip#1,#2,{\langle #1,#2\rangle}
\def\Ip#1,#2,{\bigl\langle#1,#2\bigr\rangle}
\def\IP#1,#2,{\Bigl\langle#1,#2\Bigr\rangle}

\def\mid{\,:\,}

\def\abs#1{\lvert#1\rvert}
\def\Abs#1{\bigl\lvert#1\bigr\rvert}
\def\ABs#1{\Bigl\lvert#1\Bigr\rvert}

\def\ind#1{{\mathbf 1}_{#1 } }

\def\XXint#1#2#3{{\setbox0=\hbox{$#1{#2#3}{\int}$}
     \vcenter{\hbox{$#2#3$}}\kern-.5\wd0}}

\def\md#1#2@{\ifx0#1
\begin{equation*} #2 \end{equation*}\fi  
\ifx1#1\begin{equation}#2\end{equation}\fi   
\ifx2#1\begin{align*}#2\end{align*}\fi   
\ifx3#1\begin{align}#2\end{align}\fi    
\ifx4#1\begin{gather*}#2\end{gather*}\fi  
\ifx5#1\begin{gather}#2\end{gather}\fi   
\ifx6#1\begin{multline*}#2\end{multline*}\fi  
\ifx7#1\begin{multline}#2\end{mutline}\fi  
}

\def\eqdef{\stackrel{\mathrm{def}}{{}={}}}
 
 \def\Enl #1,{\operatorname{ Enl}(#1)}
 \def\emb #1.{\operatorname{ emb}(#1)}
\def\sh#1{\operatorname{sh}(#1) }
 
\begin{document}

\title[Multiparameter Paraproducts] {Paraproducts in One and Several Parameters  }

\author{Michael Lacey}

\address{Michael Lacey\\
School of Mathematics\\
Georgia Institute of Technology\\
Atlanta,  GA 30332 USA}

\email{lacey@math.gatech.edu}

\thanks{Research supported in part by an NSF grant and a fellowship from the Guggenheim 
Foundation. The hospitality of the University of British Columbia is gratefully acknowledged.} 

\author{Jason Metcalfe}

\address{Jason Metcalfe\\
School of Mathematics\\
Georgia Institute of Technology\\
Atlanta,  GA 30332 USA}

\email{metcalfe@math.gatech.edu}

\thanks{Research supported by a NSF VIGRE grant to the Georgia Institute of Technology. 
The author is a NSF Postdoctoral Fellow. }

\begin{abstract}
For  multiparameter bilinear paraproduct operators $\operatorname B$ we prove the 
estimate 
\begin{equation*}
\operatorname B\mid L^p\times L^q \mapsto L^r, \qquad 1<p,q\le{}\infty. 
\end{equation*}
Here, $1/p+1/q=1/r$ and special attention is paid to the case of $0<r<1$.  
(Note that the families of multiparameter paraproducts are much richer than 
in the one parameter case.) 
These estimates are the essential step in the version of the multiparameter 
Coifman-Meyer theorem proved by C.~Muscalu, J.~Pipher, T.~Tao, and C.~Thiele 
\cites{camil1,camil2}.  We offer a different proof of these inequalities. 
\end{abstract}


\maketitle

\section{Introduction}  

Our subject concerns the Coifman-Meyer theorem in a multiparameter setting.  Namely, for bounded 
 function $\tau\mid \mathbb R^d \to \mathbb C$, we set 
\begin{equation*}
\operatorname T (f_1,\dotsc,f_d)(x) {}\eqdef{}\int _{\mathbb R^d} \tau(\xi) 
 \operatorname e ^{i (\xi_1+\cdots+\xi_d) x }
\prod _{j=1} ^n
\widehat f_j(\xi_j) \; d \xi_1\cdots d \xi_d
\end{equation*}
in which $f_j$ are Schwartz functions on $\mathbb R$ and $\widehat g $ denotes the 
Fourier transform, namely 
\begin{equation*}
\widehat g(\theta) {}\eqdef{}\int g(x)\operatorname e ^{-i \theta  x} \; dx.  
\end{equation*}

One is interested in  conditions under which $\operatorname T$ extends to 
a bounded multilinear operator on  a product of $L^p$ spaces. 
 And the motivation for this paper is the Theorem 

\begin{theorem}\label{t.coifmanmeyer}  Suppose that  $\tau$ obeys the 
estimates 
\begin{equation} \label{e.productsingularity}
\ABs{ \frac {\partial^\alpha}{\partial \xi^\alpha} \tau(\xi) } {}\lesssim{} 
\abs{ \xi} ^{-\alpha}=\prod _{j=1}^d \abs{ \xi_j} ^{-\alpha_j},\qquad \abs{ \alpha}\le{}N,
\end{equation}
where $\xi=(\xi_1,\dotsc,\xi_d)$ and likewise for the multiindices $\alpha$.  
There is a finite choice of $N$ so that for all $1<p_j\le{} \infty$, with not all $p_j$ being 
equal to infinity, the operator 
$\operatorname T$ extends to a bounded linear operator 
\begin{equation*}
\operatorname T\mid L^{p_1}\times\cdots\times L^{p_d}\mapsto L^r
\end{equation*}
where $1/r=\sum _{j=1}^d 1/p_j$.
\end{theorem}

In the statement of the Theorem, $A {}\lesssim{}B$ means that $A\le{}KB$ for some unspecified 
constant $K$.

This theorem has been the subject of wide ranging investigations since the inital 
results of Coifman and Meyer \cites{MR518170,MR511821}.  The methods and techniques 
of the proof, built around the subject of paraproducts, is the main focus of this article. 

The singularities permitted in $\tau$ in 
(\ref{e.productsingularity})
invoke elements of the product theory of maximal functions, singular integrals, and 
related subjects.  Some cases of the theorem above were found by Journ\'e \cite{MR88d:42028}, following the 
identification of product $BMO$ by S.--Y.~Chang and R.~Fefferman \cites{MR86g:42038,
MR82a:32009,MR90e:42030}.  

The possibility that the image space $L^r$ can have index less than one is the primary new 
contribution of C.~Muscalu, J.~Pipher, T.~Tao, and C.~Thiele \cites{camil1,camil2}.  
The purpose of this article is to give a somewhat different proof, one that discusses 
end point issues, and is a little more leisurely than the cited articles. 

The method of proving these inequalities is by way of paraproducts.  And we take the 
the latter as the primary focus of this article. See the next section for a definition 
of the most familiar paraproducts. 
Many proofs of paraproduct results depend upon the Calder\'on Zygmund decomposition, which 
has only weak analogs in the product theory. 

A very nice feature of the work of Muscalu et~al.~is that they find that the proof of 
the theorem can be understood in terms  that avoid the intricacies of the product $BMO$ theory 
of Chang and Fefferman.  We find that some aspects of that theory enter into different 
endpoint estimates, such as at $p=1$, where $L^1$ should be replaced by the Hardy space $\operatorname H^1$, 
and $p=\infty$, where the $BMO$ space enters in. 

The main result is Theorem~\ref{t.2} below, 
a discrete form of the Theorem above.  This theorem is  contained in \cites{camil1,camil2}, 
and our proof borrows elements of theirs.  We offer the proof as it differs in some 
details.  In addition, the rich family of multiparameter paraproducts is not necessarily 
well understood.   We hope that this paper advertises \cites{camil1,camil2} in particular, 
and, more generally,  the subject of multiparameter paraproducts.\footnote{A substantial part 
of the difficulties in \cites{MR1961195,witherin} is attributable to the variety of paraproducts in the multiparameter 
setting.}

\bigskip 

In the next section, we discuss the one parameter paraproducts.  The proofs in this case 
feature initial details that can be used in the multiparameter case.  We present proofs 
of these results in the special case of paraproducts formed from Haar functions. 
The subject of mutliparameter paraproducts is taken up in Section~\ref{s.multi}. 
Our presentation and proofs in Section~\ref{s.1proofs} have been influenced by the
CBMS lectures of C.~Thiele \cite{thiele}.  

 \section{One Parameter Paraproducts} 
 
 For an interval $I $, we say that $\varphi $ is \emph {adapted to $I $} iff  $\norm \varphi.2.=1$ and 
 \begin{equation} \label{e.adapted} 
 \Abs{\operatorname D ^{n} \varphi(x)} {}\lesssim{}\abs{I}^{-n-\frac12} \Bigl( 1+\tfrac {\abs{x-c(I)}}{\abs{I}}\Bigr)^{-N}, 
 \qquad n=0,1. 
 \end{equation}
 Here, $c(I) $ denotes the center of $I $, and $N $ is a large  integer,
 whose exact value  need not concern us, except to say that its value can depend upon the 
 $L^p$ inequalities that we are considering.%
\footnote{It will be clear in the sequel that $N=\min(3{p_1}+4,3{p_2}+4)$ is
 sufficient for the bilinear case, for example. 
 The main size requirement can be found in (\ref{e.adaptedN}).}
  $\operatorname D$ denotes the derivative 
 operator. 
 We shall consistently work with functions which  have  $L^2$ norm at most one.  Some of these functions we will 
also insist to have integral zero.
(Terminology for this will be introduced below.)

  Intervals will most typically be dyadic,  and we use the notation $\mathcal D$ for these intervals.  To be specific, 
 \begin{equation*}
 \mathcal D\eqdef{}\{ [j2 ^{k},(j+1)2 ^{k})\mid j,k\in \mathbb Z\}. 
 \end{equation*}
 With the control on the function and its derivative  in the definition of adapted, elements of Littlewood-Paley 
theory will apply.  Namely, we will have the inequalities (\ref{e.paleyp})---(\ref{e.paleyinfty}) 
 for the square function constructed from the functions $\{\varphi_I\mid I\in \mathcal D\} $.

 Operators are built up from rank one operators $ f\mapsto \ip f,\varphi, \varphi' $.  A paraproduct is, in its simplest manifestation,
 of the form 
 \begin{equation*}
\operatorname B(f _{1},f _{2})\eqdef{}\sum _{I\in\mathcal D} \abs{I} ^{-\frac12}\varphi _{3,I}\prod _{j=1}^{2} 
\ip f _{j},\varphi _{j,I},.
 \end{equation*}
 Here, the functions $\varphi _{j,I}$, for $j=1,2,3 $ are adapted to $I $.  Two of these three functions are 
 assumed to be of integral zero.   We should emphasize that each individual summand is of the form 
 \begin{equation*}
 (f_1,f_2)\mapsto \abs{I} ^{-\frac12}\varphi _{3,I}\prod _{j=1}^{2} \ip f _{j},\varphi _{j,I} ,. 
 \end{equation*}
 This is certainly a bounded operator from, say, $L^2\times L^2 \longrightarrow L^1 $, and our desired conclusion is that 
 the same is true for the sums above.  
 
 We will also consider higher linearities 
 \begin{equation}  \label{e.B-n-linear} 
\operatorname B(f _{1},f _{2},\ldots, f _{n})\eqdef{}\sum _{I\in\mathcal D}
\abs{I} ^{-\frac {n-1} 2}\varphi _{n+1,I}\prod _{j=1}^{n} \ip f _{j},\varphi _{j,I},
 \end{equation}
 where the functions $\varphi _{j,I}$ are assumed to be adapted to $I $ and two are of integral zero. 
 In the course of the proofs, it is convenient to consider the $n+1 $ sublinear forms 
 \begin{equation} \label{e.linearforms} 
 \operatorname \Lambda(f _{1},f _{2},\ldots, f _{n+1})\eqdef{}\sum _{I\in\mathcal D}
 \abs{I} ^{-\frac {n-1}2}\prod _{j=1}^{n+1} \abs{\ip f _{j},\varphi _{j,I},}. 
 \end{equation}
 Notice that this just assigns a number to the $n+1 $ tuple of functions and that it dominates 
 $\langle \operatorname B(f _{1},f _{2},\ldots, f _{n}), f _{n+1}\rangle $. It is also of interest to consider the related 
 sublinear operator
 \begin{equation}  \label{e.sublinearforms}
 {\operatorname L}(f _{1},f _{2},\ldots, f _{n}, f _{n+1 })\eqdef{}\sum _{I\in\mathcal D} \abs{I} ^{-\frac {n+1}2}\Biggl[\prod _{j=1}^{n+1} 
 \abs {\ip f _{j},\varphi _{j,I} , }\Biggr]\ind I .
 \end{equation}
 In particular, if $\operatorname L $ maps a product of  Banach spaces 
 into $L^1 $, then we conclude that $ \operatorname \Lambda $ is bounded on a related product of spaces, see (\ref{e.1-dual}). 
 
 \begin{theorem}\label{t.one}  For $n\ge3 $ and $1< p _{j}\le \infty $, define $\frac1r=\sum _{j=1}^{n-1}\frac1{p _{j}} $.  Then, 
 \begin{equation}  \label{e.B-1} 
 \operatorname B\mid \otimes _{j=1}^{n-1} L^{p _{j}}\longrightarrow L^r .
 \end{equation}
 In addition we have the endpoint estimates: 
 \begin{gather}   \label{e.1-1-endpoint}  
 \operatorname B\mid \otimes _{j=1}^{n-1}  L ^{1}\longrightarrow L^ {\frac 1{n-1}, \infty }
 \\   \label{e.1-bmo-endpoint} 
 \operatorname B\mid \otimes _{j=1}^{s}    L ^{1} \times \otimes _{j=1} ^{n-1-s} BMO \longrightarrow L^ {\frac 1{s}, \infty }
 \end{gather}
In this last display, we require that the functions $ \varphi_{j,I} $ have integral zero for any choice of $j $ for which $f_j $
is only assumed to be $BMO $.
 \end{theorem} 
 
 The estimates above follow immediately from the corresponding estimates for the sublinear operator, 
 in the case that the index of the 
range is between $1$ and $\infty$.  Namely for $1<r<\infty$, this can be derived from 
 
 \begin{theorem}\label{t.one-sublinear}  
 For $n\ge 3 $ and $1< p _{j}\le \infty $, define $\frac1s=\sum _{j=1}^{n}\frac1{p _{j}} $.   If 
it is the case that $0<s<\infty$, then  
 \begin{equation} \label{e.L-1}
 \norm \operatorname L(f _{1},f _{2},\ldots, f _{n} ).s. {}\lesssim{} \prod _{j=1}^{n}\norm f.p _{j}.. \end{equation}
  In the case that any $p _{j }=1 $, then $L ^{1 }$ can be 
 replaced by $H^1 $, and the estimate above is true.   If we do not replace $L^1 $ by $H^1 $, then only the 
 weak type inequality is true.  In particular, we have the estimate 
 \begin{equation} \label{e.1-L-1} 
 \operatorname L\mid \otimes _{j=1} ^{n} L^1 \longrightarrow L ^{1/n,\infty}. 
 \end{equation}
 In the case that any of $p _{j} $ equal $\infty $  and the functions $\varphi _{j,I } $  have integral zero for all $I $, 
 then the space $L ^{\infty}$ above can be replaced by $BMO $.   
 \end{theorem}

 The essential case is that of $n=3 $ above, and to avoid unnecessary notations, that is the case discussed in the proof. 
 Thus, we have three functions $\varphi _{j,I}$.  Two of these are assumed to be of integral zero.  Due to the symmetry of the 
 estimates we are to prove, we can assume that these two functions occur for $j=2,3 $.  
 
 This in particular means that we have the estimate 
 \begin{equation}  \label{e.1-maxfun}
\sup _{I\in\mathcal D }\ind I (x) \frac{\abs{\ip f_1,\varphi _{1,I},}}{\sqrt{\abs I}} {}\lesssim{} \operatorname M (f _{1})(x)
 \end{equation}
 where $M $ denotes the maximal function.   There is another bound that applies to the second and third functions.  Namely, we 
 set 
 \begin{equation}  \label{e.1-sqfn}
 \operatorname{S} _{j}g\eqdef{}{} \Bigl[ \sum _{I\in\mathcal D}\frac{\abs { \ip g,\varphi _{j,I}, } ^{2} } {{\abs I}} \ind I \Bigr] ^{1/2}.
 \end{equation}
 It is a consequence of the integral zero assumption placed on the functions $\varphi _{2,I}$ and $\varphi _{3,I}$ that the  
 usual Littlewood-Paley theory applies to these square functions. 
 Therefore, they map all $L ^{p} $ into themselves, for $1<p< \infty $, and we have the usual endpoint estimates. 
 To be explicit, these estimates are 
 \begin{align}  \label{e.paleyp}
 \norm \operatorname S_j f.p. &{}\lesssim{}\norm f.p.,\qquad 1<p<\infty,
 \\ \label{e.paleyH1}
 \norm \operatorname S_jf.1. &{}\lesssim{}\norm f. H^1.,
 \\ \label{e.paleyinfty}
 \norm \operatorname S_jf.BMO. &{}\lesssim{}\norm f.\infty. .
 \end{align}
 
 Note that as we are using the maximal function and square functions, we have access to the following upper bound for 
 numerical sequence 
 \begin{equation} \label{e.zI-2-2}
 \sum _{n} \prod _{j=1} ^{3} a _{j,n}\le{ } \norm a _{1,n}. \ell^ \infty.
 \norm a _{2,n}.\ell^2. \norm a _{3,n} .\ell^2..
 \end{equation}
 
 \subsection{Generalities on the Proof} 
 
  If $1<{} p_j < \infty $ and $\frac 1{p_1}+\frac 1{p_2 }+\frac1{p _{3 } }=1$, we can estimate, based on (\ref{e.zI-2-2}) 
  and H\"older's  inequality,  
  \begin{equation} \label{e.1-dual} 
  \begin{split} 
 \int \operatorname B(f_1,f_2) \overline{f_3}\; dx & \le{} 
 \int \operatorname L(f_1,f_2,f_3)\; dx
 \\&{}\le{} \int(\operatorname{M}f_1)\cdot(\operatorname{S}_2f_2)\cdot(\operatorname{S}_3f_3)\; dx
  \\&{}\le{} 
  	\norm \operatorname{M}f_1.p_1.\norm \operatorname{S}_2f_2.p_2.\norm \operatorname{S}_3f_3.p_3.
\\& {}\lesssim{} 
  	\norm  f_1.p_1.\norm  f_2.p_2.\norm  f_3.p_3..
	\end{split} 
	\end{equation}
 This argument also applies when $p_1=\infty$.   
 When, however,   any of the 
 $p_j=1 $, one should  replace the $L^1 $  norm on $f_j$ with the  $H^1 $ norm.  
 
 The argument must be modified when, e.g.,~$p_2 =\infty $. For then 
 the square function $\norm \operatorname{S}_2f_2.\infty. $ is no longer bounded.  And indeed, 
 the sharp estimate on the square function replaces $L^\infty $ with $BMO $.

  Alternate methods are required when duality cannot be applied.  Here, we shall 
  obtain inequalities of weak type.  For example, 
  \begin{equation}  \label{e.1-WEAK-1}
  \lambda\abs { \{ \operatorname B(f_1,f_2)>\lambda \} } ^{1/r} {}\lesssim{}\norm f_1.p_1.\norm f_2.p_2.,\qquad 
  	\tfrac 1{p_1}+\tfrac 1{p_2 }=\tfrac1r.
  \end{equation}
  Interpolation will then supply the strong type inequalities, except for the 
  endpoint estimates.
  
  As the class of operators $\operatorname B $ we consider are invariant under   dilations by powers of $2 
  $, this inequality follows from 
  \begin{equation}  \label{e.1-weak-1}
  \abs{\{\operatorname B(f_1,f_2)>K\} } {}\le{} 1
  \end{equation}
  where $K $ is an absolute constant, and the inequality holds for all choices of smooth compactly supported functions  $f_j $
  with $L ^{p_j }$ norm $1 $.

  The usefulness of this observation is already evident in that we have the following (obvious) estimate 
  \begin{equation}  \label{e.weak-exercise}
  \norm fg.1/2,\infty. {}\lesssim{}\norm f.1,\infty.\norm g.1,\infty..
  \end{equation}
  This inequality immediately generalizes to general products and indices.  We use this generalization below.
  
  One can effectively use the symmetry in the formulation of the 
  paraproducts in passing to the sublinear function $\operatorname \Lambda$, 
  and considering weak type inequalities for it. 
     Namely, for $1\le{} p_1,p_2\le{}\infty$, we 
  define $p_3 $ by $\frac 1{p_1}+\frac 1{p_2 }+\frac1{p _{3 } }=1$.  In particular, $p_3 $ can be negative:  For $p_1=p_2=1 $, 
  we have $\frac1{p_3}=-1 $, which we interpret as the dual index to $\frac12 $.  Let $X(E) $ be the space of functions supported on a measurable set $E\subset \mathbb R $ and bounded by 
  $1 $.  
  We then prove the inequality 
  \begin{equation}  \label{e.1-LambdA-1}
  \abs{\operatorname \Lambda(f_1,f_2,f_3)} {}\lesssim{}\abs{ E_3} ^{1/p_3}\prod _{j=1}^2 \norm f_j. p_j. 
  \qquad f_3\in X(E_3). 
  \end{equation}
  Observe that this implies (\ref{e.1-weak-1}).  Also observe that the inequality 
  for $\operatorname \Lambda$ follows from the 
   following  formulation:  
  For all $E_3 $, we can choose $E_3'\subset E_3 $ with  $\abs{E'_3}\ge\frac12\abs{E_3} $, and 
  \begin{equation}  \label{e.1-Lambda-1}
  \abs{ \operatorname \Lambda(f_1,f_2,f_3)} {}\lesssim{}
  \abs{ E_3} ^{1/p_3}\prod _{j=1}^2 \norm f_j. p_j. 
  \qquad f_3\in X(E_3'). 
  \end{equation}
  By dilation invariance, it suffices to prove this estimate in the case that 
  $\norm f_1.p_1.=\norm f_2.p_2.=1$ and $\abs{ E_3}=1$.
  All of these comments apply equally well in the multiparameter case. 
 
 \subsection{ $H^1$ and $BMO $}
 
 We will restrict ourselves to the dyadic versions of the real Hardy space $H^1 $ and its dual $BMO $.

 The Haar functions are 
 \begin{equation}  \label{e.haar}
 h_I = \abs{ I} ^{-1/2}(\mathbf 1 _{I_-}-\mathbf 1 _{I_+}),\qquad I\in \mathcal D,
 \end{equation}
 where $I_-$ ($I_+$) is the left (right) half of $I$.  These functions form a basis for $L^2$. 
 The dyadic square function from Haar functions is formed as follows. 
 \begin{equation*}
 \operatorname Sf {}\eqdef{}\Bigl[\sum_{I\in \mathcal D} \frac{\abs{ \ip f,h_I,}^2} {\abs{ I} } \mathbf 1 _{I} \Bigr] ^{1/2}.
 \end{equation*}

 We define the real dyadic Hardy space $H^1$ as those functions $f$ with 
 \begin{equation*}
 \norm f .H^1. {}\eqdef{}\norm f.1.+\norm \operatorname S f.1.<\infty.
 \end{equation*}
 
 The dual to $H^1$ is $BMO$.  This space has the equivalent norm 
 \begin{equation}  \label{e.bmo}
 \norm f.BMO. {}\eqdef{} \sup _{J \in \mathcal D} \Bigl[ \abs{ J} ^{-1} \sum_{ I\subset J  } \abs{ \ip f,h_I,}^2\Bigr] ^{1/2}.
 \end{equation}
 
 These  spaces are substitutes for $L^1$ and $L^\infty$.  In the current setting sharp endpoint estimates can 
 be phrased in terms of these norms.  And there is a rich interpolation theory between these spaces.

 \section{Proofs in the One Parameter Case} \label{s.1proofs}

  The case of Haar paraproducts is the only case that we consider in the one parameter case. 
  The rationale is the proof in the multiparameter case includes the one parameter case 
  as a special instance.  In addition, the Haar case 
  is especially attractive, due to the presence of the dyadic grid.  
  
  A particular way that it enters is this. 
  Suppose that $\mathcal I $ is a  collection of  disjoint dyadic intervals, not necessarily a partition of $\mathbb R $.
  We define the conditional expectation   with respect to $\mathcal I $ as 
  \begin{equation*}
  \mathbb E (f \, | \,\mathcal I )(x)\eqdef{}{} 
  	\begin{cases} 
	   \abs{ I} ^{-1}\int_I f(y)\; dy  & x\in I,\  I\in\mathcal I, 
	   \\ 
	   f(x)		& x\not\in\bigcup _{I\in\mathcal I }I
	   \end{cases}
  \end{equation*}
  We leave it as an exercise that these properties of the conditional expectation are true. 
  \begin{enumerate} 
  \item  { Integrals are preserved under conditional expectation: $\int f \; dx=\int \mathbb E(f\, | \, \mathcal I ) \; dx $.} 
  \item{ $f \mapsto\mathbb E (f \, | \,\mathcal I ) $ is a projection.} 
  \item { $f \mapsto\mathbb E (f \, | \,\mathcal I ) $ is of norm one on all $L^p $, $1\le{}p\le{}\infty $. }
  \item{$f \mapsto\mathbb E (f \, | \,\mathcal I ) $ is bounded as a map from dyadic $H^1 $ into itself. }
  \end{enumerate} 
  
  \bigskip 
  
  We first turn to the range of inequalities for the sublinear operator $\operatorname L $ 
and the proof of Theorem~\ref{t.one-sublinear}. 
  Observe that by (\ref{e.zI-2-2}), we have 
  \begin{equation}  \label{e.trivial-on-L}
   \operatorname L(f_1,f_2,f_3)\le{}\operatorname M f_1\cdot \operatorname S f_2\cdot \operatorname S f_3 
  \end{equation}
  Here, we assume that we have mean zero in the second and third places, and we continue with this assumption below. 
To be specific, the sublinear operator is 
\begin{equation*}
\operatorname L(f_1,f_2,f_3)=\sum _{I\in \mathcal D} \mathbf |I|^{-3/2} \mathbf 1 _{I}\abs{ \ip f_1,\abs{h_I}, 
\ip f_2,h_I, \ip f_3,h_I, }.
\end{equation*}
For $f_1$, we form the inner product with the absolute value of the Haar function. 
  The inequalities in (\ref{e.L-1}) then follow from H\"older's 
inequality, provided that all we are not discussing an endpoint estimate.  When $s<1$, one can instead apply
an appropriate version of (\ref{e.weak-exercise}).

If any $p_j=1 $, then we only conclude that $ \operatorname M f_j $ and $\operatorname S_j f_j $ are in $L^{1,\infty} $. 
But we can apply (\ref{e.weak-exercise}) to conclude the weak type estimate.  
If any $p_j=1 $ and $f_j\in H^1 $, then we conclude that both $\operatorname M f_j $ and $\operatorname S_j f_j $ 
are in $L^1 $, so that again H\"older's inequality or (\ref{e.weak-exercise}) will apply.  

We concern ourselves with the  endpoint estimates  where 
either of $p_2,p_3 $ is infinity and $L^\infty $ is replaced with $BMO $.  
One class of inequalities are in fact easily available; they are 
\begin{equation} \label{e.x-1-estimate}
		\begin{split} 
\operatorname L &\mid L^\infty \otimes BMO \otimes H^1\longrightarrow L^1,
\\
\operatorname L &\mid L^\infty \otimes BMO \otimes BMO\longrightarrow BMO.
		\end{split}
\end{equation}
Notice that these estimates can be interpolated by standard linear methods.    

Since the Haar functions are an unconditional basis for both $H^1 $ and $BMO$,\footnote {While we are specifically 
appealing to the properties of the Haar functions here, this aspect does generalize to the non-Haar functions. }
we can conclude that 
\begin{align*}
\int \operatorname L (f_1,f_2,f_3)\; dx& {}={}
\sum _{I\in\mathcal D } \frac {\abs {\ip  f_1, |h_I| , }} {\sqrt {\abs I } } 
\abs {\ip f_2, h_I , } \abs {\ip f_3, h_I , } 
\\&{}\lesssim{} 
\norm f_1.\infty.\norm f_2 . BMO . \norm f_3 . H^1 .. 
\end{align*}
This proves the first bound. 

For the $BMO $ estimate, for each dyadic interval $J $ we have 
\begin{align*}
\sum _{I\subset J } \frac {\abs {\ip  f_1, |h_I| , } } {\sqrt {\abs I } } \abs {\ip f_2, h_I , } \abs {\ip f_3, h_I , } &{}\lesssim{} 
\norm f_1.\infty.   \prod _{j=1}^2  \Bigl[ \sum _{I\subset J } \abs {\ip f_2, h_I , }^2 \Bigr] ^{1/2 }
\\&{}\lesssim{} \abs  J \,\norm f_1.\infty. \prod _{j=1} ^2 \norm f_j . BMO . . 
\end{align*}
This concludes the proof of the estimates (\ref{e.x-1-estimate}).

\medskip

The last estimates to prove are these: 
\begin{equation*}
\operatorname L \mid L ^{p_1 }\otimes BMO \otimes L ^{p_3 }\longrightarrow L ^{s },\qquad 
\tfrac1s=\tfrac1{p_1 }+\tfrac1{p_3 }.
\end{equation*}
At this point we make a more substantive reliance on the dyadic structure.   The strategy is first to  
prove the weak type inequalities, and in particular (\ref{e.1-LambdA-1}).  Namely, we will choose an 
exceptional set on which we will not attempt 
to estimate $\operatorname L (f_1,f_2,f_3) $ and a conditional expectation to apply to $f_1 $, after which we will have a 
bounded function in the first coordinate.  But then we will be in a situation for which we can appeal to 
the estimates in which we have duality.

We will prove that $\operatorname L $ satisfies (\ref{e.1-LambdA-1}).  Thus, fix $f_j $ functions in the appropriate spaces, of norm 
one.  Define 
\begin{equation*}
E\eqdef{}\{ \operatorname M f_1>1 \}.
\end{equation*}
We do not attempt to estimate $\operatorname L $ on this set.  That has the practical implication that we need only 
consider the sum 
\begin{equation*}
\operatorname L _E (f_1)\eqdef{}\sum_{I\not\subset E} |I|^{-1}
\frac {\abs{ \ip f_1, \abs{ h_I} , }} {\sqrt {\abs I } } \abs {\ip f_2, h_I , } \abs {\ip f_3, h_I , } \ind I 
\end{equation*}
(Recall that we are assuming that $f_2 $ and $f_3 $ are fixed.) 
Let $\mathcal I $ be the collection of maximal dyadic subintervals of $E $, and set $g_1\eqdef{}\mathbb E(f_1\, | \,\mathcal I) $.

Note that by construction we will have the estimate 
$\norm g_1.\infty.\le{}2 $. For otherwise, let $J$ be the smallest dyadic interval that strictly contains 
$I$ (i.e.~the parent of $I$), and observe that $\int _{J}\abs{g}\; dx\ge\abs J$.  That is, we contradict the 
maximality of $I$.

In addition, for each dyadic interval $I $ not contained in $E $, we have 
$ \int_I f_1\; dx=\int _I g \; dx $.  Thus, we have 
\begin{equation}  \label{e.conditional-trick}
\operatorname L _E (f_1)=\operatorname L _E (g_1).
\end{equation}
Therefore, using (\ref{e.x-1-estimate}), we can estimate 
\begin{align*}
\abs{\{ \operatorname L (f_1,f_2,f_3)>1\} }&{}\le{} 
			\abs E+\abs{ \{\operatorname L (g_1,f_2,f_3)>1\}}
\\&{} {}\lesssim{} {} 
		1+{} [\norm g_1.\infty. \norm f_2.BMO. \norm f_3.p_3.]^{p_3}
\\& {}\lesssim{} 1.
\end{align*}
Our discussion of the estimates in Theorem~\ref{t.one-sublinear} is complete.

\bigskip 

Let us turn to the bilinear operator $\operatorname B (f_1,f_2) $ and the proof of Theorem~\ref{t.one}.  In 
the inequality (\ref{e.B-1}), if the index $r $ of the target space  is  between $1 $ and $\infty $, then 
we can appeal to 
duality, as is done explicitly in (\ref{e.1-dual}).

We discuss the proof of the weak type bounds for $\operatorname B $, in the case that duality does not apply, namely
$\frac 12\le r<1 $. 
Marcinkiewicz interpolation will then deduce the strong type $L^r $ inequalities. 

In so doing, we need only prove (\ref{e.1-weak-1}), and we will repeat the use of conditional expectation 
in the argument (\ref{e.conditional-trick}) above. 
 Take $f_j \in L^{p_j } $ of norm one, for $j=1,2 $.   
 Suppose that we have in fact $\norm f_j .\infty.\le1 $.  
 We conclude that in fact $\norm f_j.q.\le1 $ for all $p_j<q<\infty $, and so for $ q >4$ large, we can use the proven 
 bound of $L^q \times  L^q $ into $L ^{q/2} $ to conclude that 
 \begin{equation*}
\abs{  \{ \operatorname B(f_1,f_2)>K\} } {}\le{}1. 
 \end{equation*}
   The general case can be  reduced to this situation.

   Define 
 \begin{equation*}
 E\eqdef{}\bigcup _{j=1}^2 \{ \operatorname M f_j>1\},\qquad F\eqdef{}\{\operatorname M \ind E>\tfrac12\}.
 \end{equation*}
 Clearly, the set $F $ has measure bounded by an absolute constant. 
 We will not estimate $\operatorname B (f_1,f_2) $ on the set $F $.   Define 
 \begin{equation*}
 \operatorname B _F(f_1,f_2)=\sum _{I\not \subset F }\abs I ^{-1/2 } h_I 
 \ip f_1 ,\abs{h_I} ,\ip f_2 ,h_I ,. 
 \end{equation*}
 And set $\mathcal I $ to the collection of maximal dyadic intervals contained in $F $.  
We set $g_j\eqdef{}\mathbb E(f_j\, | \,\mathcal I) $. 
 Then, certainly we have $\norm g_j .\infty.\le1 $.    We claim that 
 \begin{equation} \label{e.1-=} 
 \operatorname B _F(g_1,g_2) \ind {F^c}=\operatorname B(f_1,f_2)\ind {F^c}. 
 \end{equation}
 And this will complete our proof.   
 
 Suppose that $I $ is a dyadic interval that is not contained in $F $.  The Haar function associated to $I $ is constant 
 on the two sub halves of $I $, which we denote as $I_\pm $.  By our definition of $F $, neither $I_\pm $ can be contained in $E $, 
 hence we have 
 \begin{equation*}
 \int_{I_\pm} f_j\; dx=\int _{I_\pm } g_j \; dx. 
 \end{equation*}
 This proves (\ref{e.1-=}), and so we have completed the proof of the norm bounds for $\operatorname B $.

   \section{Multiparameter Paraproducts}  \label{s.multi}
 
We now consider paraproducts formed over sums of dyadic rectangles in  $\mathbb R^d$.
The class of paraproducts is then 
invariant under a $d$ parameter family of dilations, 
a situation that we refer to as ``multi--parameter.''\footnote{In this paper, the number of dimensions 
will be the number of parameters.  In general, the two are however distinct. Consider $\mathbb R ^{d_1}\otimes 
\mathbb R ^{d_2}$ and rectangles in this space formed from a cube in each space $\mathbb R ^{d_j}$. }

Let us say that a function $\varphi $  \emph { is adapted to a rectangle $R= \mathop\otimes _{j=1}^dR _{j}$ } iff 
$ \varphi (x_1,\dots ,x_d)=\prod _{j=1}^d\varphi_j(x_j) $, with each  $\varphi_j $ adapted to the interval $R_j $ 
in the sense of (\ref{e.adapted}).

Our paraproducts are of the same general form
\begin{equation*}
\operatorname B(f _{1},f _{2},\ldots, f _{n})\eqdef{}\sum _{R\in\mathcal R}
\frac{\varphi _{n+1,R}}{\abs{R}  ^{\frac{n-1}2}}\prod _{v=1}^{n} 
\ip f _{v},\varphi _{v,R},.
 \end{equation*}
 Here, we let $\mathcal R \eqdef{}\mathcal D^d $ be the class of dyadic rectangles. 
 With the obvious changes, we will 
 also use the notations for the sublinear forms and operators given in (\ref{e.linearforms}) and (\ref{e.sublinearforms}). 
 
 The Theorem in this setting is 
 
 \begin{theorem}\label{t.2}  Let $n\ge2 $ and  $1<p _{v }\le{} \infty $ for 
 $1\le v\le{} n $, and define $\frac1r=\sum _{v=1 }^{n }\frac 1{p _{v } }$. 
 Assume that for each choice of coordinate $1\le j \le d$,  
 there are two choices of $1\le v \le n+1$ for which we have 
 \begin{equation}  \label{e.2-zeros}
 \int _{\mathbb R }\varphi _{v,R }(x_1,x_2,\cdots,x_n)\; dx_j=0,\qquad \text{for all $x_ k$ with $k\not=j$ and all $R $ }.  
 \end{equation}
 Then, we have the inequality 
 \begin{equation}  \label{e.B-2-}
 \operatorname B\mid \otimes _{v=1}^{n} L^{p _{v}}\longrightarrow L^r .
 \end{equation}
 Assume that 
 the functions $ \varphi _{v,R}$ 
  satisfy (\ref{e.2-zeros}) for all $ j$.  In the instance that 
 $ p_v=1$, the inequality remains true if we replace $ L^1$ by $ H^1$ defined below. 
 In the instance that $ p_v=\infty $, we can replace $ L^\infty $ by the larger 
 space $ {BMO}=(H^1) ^{\ast}$ defined below. 
 \end{theorem}

  The critical distinction comes from the assumption about the zeros, (\ref{e.2-zeros}).  
 Let us say that \emph { there are  $x _{j } $ zeros in the $v $th position } iff 
\begin{equation}  \label{e.x1-zeros}
\int \varphi _{v,R }(x_1,\dotsc ,x_d)\; dx _{j }=0\qquad\text{for all $x _{k }$ with $k\not=j$. }
\end{equation} 
 And so our assumption is that for each $1\le j\le d$ there are two choices of $v$ for 
 which we have zeros in the $v$th position.

 One would expect that in the Theorem above the space $L^\infty$ could be replaced by the 
 larger product $BMO$ space.  This is a more delicate consideration, one that we remark upon 
 at the end of the paper.

 Again, the critical case is $n=2 $, so that $\operatorname B$ is bilinear. 
 There are essentially $d$ distinct cases.  The first case, with the greatest similarity to the one parameter case, is where 
 we have, for example,  $x_j $ zeros in first and second positions for all $1\le j\le d$.   
 The other cases do not have a proper analog in the one parameter case. 
 \subsection{ $H ^{1 }$ and $BMO $} 
  We turn to the product Hardy space theory, as developed by S.-Y.~Chang and R.~Fefferman \cites{MR86g:42038
,
MR82a:32009}.   This section is not strictly speaking needed, but does inform 
the modes of proof below.

$\operatorname H^1(\mathbb C_+^d)$ will denote  the 
$d$--fold product  real valued Hardy space.  
This space consists of functions $f\mid 
\mathbb R^d\longrightarrow\mathbb R$ where  $\mathbb R^d$ is viewed as the boundary of 
\begin{equation*}
\mathbb C_+^d=\otimes_{j=1}^d \{ z\in\mathbb C\mid \operatorname{Re}(z)>0\}.
\end{equation*}
We require that there is a function $F\,:\,\mathbb C_+^d\longrightarrow\mathbb C$ that is holomorphic in each variable separately 
and 
\begin{equation*}
f(x)=\lim_{\norm y..\to0}\operatorname {Re}F(x_1+iy_1,\ldots,x_d+iy_d).
\end{equation*}	
The norm of $f$ is taken to be 
\begin{equation*}
\norm f. \operatorname H^1 .=\lim_{y_1\downarrow 0}\cdots\lim_{y_d\downarrow 0}\norm
F(x_1+iy_1,\ldots,x_d+iy_d).L^1(\mathbb R^d )..
\end{equation*}

Product $\operatorname H^1(\mathbb C_+^d)$ has the equivalent norm 
\begin{equation} \label{e.H1def}
\norm f.\operatorname H ^1. {}\eqdef{} \norm f.1.+\norm \operatorname Sf.1.
\end{equation}
where $\operatorname S$ is the square function formed over the product Haar basis 
\begin{equation*}
\operatorname S f {}\eqdef{} \Bigl[\sum_{R\in \mathcal R} \frac{ \abs{ \ip f, h_R,}^2 }{\abs{ R}} \mathbf 1 _{R} \Bigr] ^{1/2}.
\end{equation*}
For a rectangle $R=\prod_{j=1}^d R_{(j)}\in\mathcal D^d$, we have
set 
\begin{equation*}
h_R(x_1,\ldots,x_d)=\prod_{j=1}^d h_{R_{(j)}}(x_j).
\end{equation*}
The last product is over one dimensional Haar functions as in (\ref{e.haar}). 
The basis $\{ h_R\mid R\in\mathcal D^d\}$ is the $d$--fold tensor product of the Haar basis.

The dual of  $\operatorname H^1(\mathbb C_+^d)$ is $\operatorname{BMO}(\mathbb C_+^d)$, the $d$--fold product $\operatorname{BMO}$ space. It is a Theorem of S.-Y.~Chang and R.~Fefferman 
\cite{MR82a:32009} that this space has an explicit characterization in terms of  the product Haar basis.  
In particular,
Chang and Fefferman showed that
 the product $\operatorname{BMO}$ space has the equivalent norm 
\begin{align*}
\lVert b\rVert_{\operatorname{BMO}}={}&\sup_{ U\subset\mathbb R^d } \Bigl[
\abs{{ U}}^{-1} \sum_{R
\subset
 U}\abs{\ip b,h_R,}^2\Bigr] ^{1/2}
\end{align*}
where it essential that the supremum be formed over all subsets $U\subset \mathbb R^d$ 
of finite measure. 

\subsection{The Governing Operators }
We describe a range of operators, which encompass the $d$ parameter maximal function 
at one end and the $d$ parameter square function at the other.  These operators, 
as we shall see, govern the behaviors of these paraproducts.

  To be explicit, in the two parameter setting, these operators are as follows. First 
  we have the maximal function, 
  \begin{equation*}
  \operatorname{MM} f\eqdef{}\sup _{R\in\mathcal R} \frac { \abs{\ip f,\varphi _{R }, }} {\sqrt{\abs R} } \ind R ,
  \end{equation*} 
  which is a variant of the strong maximal function.  
  
  The reason for the iterated style 
  notation becomes clearer with the second type of governing operator.  It is 
  \begin{equation*}
  \operatorname {S_1M_2} f {}\eqdef{} \Bigl[
  \sum _{R_{(1)}\in \mathcal D}\sup _{R_{(2)}\in \mathcal D} \frac{ \abs{\ip f,\varphi_{R _{(1)}\times R _{(2)}},}^2}
  {\abs{R } } \mathbf 1 _{R}\Bigr] ^{1/2}, \qquad R=R_{(1)}\times R_{(2)}.
  \end{equation*}
  In order for this to be a bounded operator, the functions $\{\varphi _{R}\}$ must 
  have zeros in the first coordinate.  But then, the operator will be bounded 
  on all $L^p$'s for $1<p<\infty$. There is also the operator $\operatorname {S_2M_1}$ 
  in which the role of the coordinates is changed. 
  
  A third type of operator is 
  \begin{equation*}
  \operatorname {M_1S_2} f {}\eqdef{} \sup _{R_{(1)}\in \mathcal D} \Bigl[
  \sum _{R_{(2)}\in \mathcal D}\frac{ \abs{\ip f,\varphi_{R _{(1)}\times R _{(2)} },}^2}
  {\abs{R } } \mathbf 1 _{R}\Bigr] ^{1/2}, \qquad R=R _{(1)}\times R _{(2)}.
  \end{equation*}
  The functions  $\{\varphi _{R}\}$ must now 
  have zeros in the second coordinate. And there is a corresponding operator $\operatorname {M_2S_1}$
  in which the roles of the coordinates are reversed. 
  
  A fourth type of operator is the familiar two parameter square function 
  \begin{equation*}
  \operatorname {SS}f {}\eqdef{}\Bigl[\sum _{R\in \mathcal R} \frac{ \abs{\ip f,
  \varphi_{R},}^2}
  {\abs{R } } \mathbf 1 _{R} \Bigr] ^{1/2}.
  \end{equation*}
  Here, we require that the functions $\varphi_R$ have zeros in both coordinates. 
 As with the maximal function $\operatorname M \operatorname M$, 
 the subscripts are not needed in this case.
  
  \smallskip 
  
In general, we set $\operatorname T_j$  to be either the square function $\operatorname S$
or the maximal function $\operatorname M$, both formed over a set of functions 
$\{\varphi_I\mid I\in \mathcal D\}$ acting on the $j$th coordinate.
For a permutation of the coordinates 
$\pi\mid \{1 ,\dotsc ,d\}\mapsto \{1 ,\dotsc ,d\}$, set 
\begin{equation}  \label{e.U}
\operatorname T {}\eqdef{}\operatorname T _{\pi(1)}\cdots \operatorname T _{\pi(d)}. 
\end{equation}
The subscript $\pi(j)$ indicates in which coordinate the operator $\operatorname T _{\pi(j)}$ 
operates.  
In each position in which coordinate $\operatorname T _ {\pi(j) }$ 
is the square function, we require that the functions 
$\{\varphi_R\mid R\in \mathcal R\}$ have zeros in that coordinate. 
Note that these operators can be viewed as 
\begin{equation} \label{e.repeatednorms}
\operatorname T f(x) =\Bigl\lVert \cdots \Bigl\lVert \Bigl\{ 
\frac  {\abs{ \ip f, \varphi _{R  } ,} } 
{\sqrt{\abs{ R}} } \mathbf 1 _{R }\mid R=R _{(1)}\times \cdots \times R _{(d)}
\in \mathcal R \Bigr\}  \Bigr\rVert _ {\ell ^{\sigma(\pi(d))}( R _{(\pi(d))}) }
\cdots \Bigr\rVert _ {\ell ^{\sigma(\pi(1))}( R _{(\pi(1))}) }
\end{equation}
where $\sigma(j) $ is either $2$ or $\infty$ for all $j$.

The point of these definitions is that for all paraproducts $\operatorname B$ of $d$ parameters,
there are three choices of $\operatorname T_k$, $k=1,2,3$, operators as in 
in (\ref{e.U}) for which we have 
\begin{equation}  \label{e.ipB}
  \ip {\operatorname B(f_1,f_2) }, f_3 ,\le{} 
  \int \prod _{k=1}^3 \operatorname T_k(f_k)\; dx.
  \end{equation}  
 This is a consequence of the essential hypothesis on there being two zeros
 in each coordinate.  In those two positions, one uses the square function. 
 In every other position, the maximal function is used.

\subsubsection*{$L^p$ bounds for the operators $\operatorname T$.}
Let us discuss the mapping properties of these operators, beginning with 
the maximal operator. We have been careful to insist that the functions 
$\varphi_R$ are products of adapted functions.  Thus, in the two parameter case, 
we can appeal to the one parameter maximal function twice, as follows. 
\begin{align*}
\norm \operatorname {M_1M_2}f.p. & {}\lesssim{} 
\norm \operatorname {M_1} \{ \operatorname M_2 f \} .p.\\
& {}\lesssim{} \norm \operatorname M_2 f .p.
\\ & {}\lesssim{}\norm f.p..
\end{align*}

Likewise, by a $d$ fold iteration of this argument, it follows that 
\begin{equation*}
\norm \operatorname {M\cdots M} f.p. {}\lesssim{}\norm f.p.,\qquad 1<p\le\infty. 
\end{equation*}

The same estimates are true for the square function, but are not as straightforward to deduce.  

\begin{lemma}\label{l.SS}  Assume that the functions $\{\varphi_R\}$ have zeros 
in every coordinate. 
 Then we have the inequalities below
\begin{align}  \label{e.SSp}
\norm \operatorname {S\cdots S} f.p. &{}\lesssim{}\norm f.p.,\qquad 1<p<\infty, 
\\   \label{e.SS1} 
\norm \operatorname {S\cdots S} f.1. &{}\lesssim{}\norm f.H^1. .
\end{align}
At the $L^\infty$ endpoint, the correct estimate is
\begin{equation}
\label{e.SSbmo} \sup _{U} \abs{ U} ^{-1}\sum _{R\subset U}
\abs{\ip f, \varphi_R ,}^2{}\lesssim{}\norm f.L^\infty.^2.
\end{equation}
The supremum is formed over all subsets  $U \subset \mathbb R^d $ of finite measure. 
\end{lemma}

\begin{proof}
We should  consider the one parameter inequality 
\begin{equation} \label{e.12}
\norm \operatorname S f. 2. {}\lesssim{}\norm f.2., 
\end{equation}
as this will explain in part the assumptions used in the definition of adapted (\ref{e.adapted}).

Consider first the inner product $\ip \varphi_I,\varphi_J,$ for two dyadic intervals $\abs{ I}\le{} \abs{ J}$. 
Using the fact that $\varphi_I$ has integral zero and that $\varphi_J$ admits a control on 
it's first derivative, we estimate 
\begin{align*}
\rho(I,J)& {}\eqdef{}\abs{ \ip \varphi_I,\varphi_J,} 
\\&{}\le{}  \int \abs{ \varphi_I(x)[\varphi_J(x)-\varphi_J(c(I))]}\; dx 
\\& {}\le{} \bigl(\tfrac {\abs{ I}} {\abs J }\bigr) ^{3/2} 
\bigl(1+\tfrac{\abs{ c(I)-c(J)} } {\abs J } \bigr) ^{-2}.
\end{align*}
Here recall that $c(I)$ is the center of $I$.  For this particular argument we only need 
$N=3$, say, though other parts of the argument require higher values.

Now, observe that we have 
\begin{align*}
\rho_1 &{}\eqdef{}\sup _{I}\sum _{J\mid \abs I\le{}\abs J} \rho(I,J)<\infty,
\\
\rho_2 &{}\eqdef{}\sup _{J}\sum _{I\mid \abs I\le{}\abs J} \rho(I,J)<\infty.
\end{align*}

To prove (\ref{e.12}), observe that by duality, it suffices to prove the estimate 
\begin{equation*}
\NOrm \sum _{I}a_I \varphi_I .2. {}\lesssim{}\Bigl[\sum _{I}\abs{ a_I}^2 \Bigr] ^{1/2}. 
\end{equation*}
Assume that the right hand side is one, and  estimate
the square of the left hand side as follows, with a generous use of the  Cauchy-Schwarz  
inequality. 
\begin{align*}
\NOrm \sum _{I}a_I \varphi_I .2. &\le{}2\sum _{\abs I\le{}\abs J}\abs{ a_I a_J } \rho(I,J) 
  \\&{}\le{}\Bigl[\sum _{I} \ABs{ \sum _{J\mid \abs I\le{}\abs J} \abs{ a_J} 
  \rho(I,J) }^2\Bigr]^{1/2}
  \\&{}\le{}2 \Bigl[ \sum _{I}   \Bigl\{ 
  \sum _{J\mid \abs I\le{}\abs J}\abs{ a_J}^2 \rho(I,J)\Bigr\} 
   \Bigl\{ \sum _{J\mid \abs I\le{}\abs J}\rho(I,J)\Bigr\}\Bigr] ^{1/2}
  \\&{}\le{} 2\rho_1 \rho_2. 
\end{align*}
This completes the proof of (\ref{e.12}).  

 \medskip

The proof of (\ref{e.SSp}) in the case of $p=2$ follows from an iteration of the 
one parameter result, just like the argument for the maximal function.

The $H^1$ to $L^1$ estimate is an easy consequence of the definition 
of the $ H^1$ norm in (\ref{e.H1def}).

The last estimate (\ref{e.SSbmo}) is not as accessible, as it relies upon a ``localization lemma''{} 
Lemma~\ref{l.local}.   Fix a set $U\subset \mathbb R^d$ of finite measure and a function $f$ bounded 
uniformly by one.    

Set $\operatorname T_0$ to be the $d$ parameter maximal function, and define a sequence of functions $f_k$ 
by taking 
\begin{equation*}
f_0 {}\eqdef{} f \mathbf 1 _{\{\operatorname T_0 \mathbf 1 _{U}>2 ^{-1} \} } .
\end{equation*}
And then inductively define $f_k$,  $k\ge1$ by 
\begin{equation} \label{e.inductive}
f_0+\cdots+f_k=f \mathbf 1 _{\{\operatorname T_0 \mathbf 1 _{U}>2 ^{-1-k} \} } .
\end{equation}

Of course we have 
\begin{equation*}
\sum _{R\subset U} \abs{ \ip f_0,\varphi_R,}^2 {}\lesssim{} \abs{ U}
\end{equation*}
by the $L^2$ boundedness of the square function.    
By Lemma~\ref{l.local}, for integers $k\ge1$ we have 
\begin{equation} \label{e.zz}
\sum _{R\subset U} \abs{ \ip f_k,\varphi_R,}^2 {}\lesssim{} 2 ^{-k}\abs{ U}
\end{equation}
for an appropriate choice of integer $N$ in the definition of adapted, (\ref{e.adapted}).  
As this last estimate is summable in $k$, the proof is complete. 

\end{proof}

  The case of general operators $T$ is treated in this proposition.

  \begin{proposition}\label{p.Tbounded} All possible operators $\operatorname T$ 
  as in (\ref{e.U}) map $L ^{p}$ into itself for all $1<p<\infty$.  This holds provided 
  the functions $\{\varphi_R\mid R\in \mathcal R\}$ have zeros in each coordinate 
  in which $\operatorname T$ is equal to a square function.  The norm depends only on the 
  constants that enter into the definition of adapted in (\ref{e.adapted}).
  \end{proposition}

  \begin{proof}
  In the two parameter case, observe that 
  \begin{equation*}
  \operatorname M_1 \operatorname S_2\le{} \operatorname S_2 \operatorname M_1,
  \end{equation*}
  the inequality holding pointwise.  
  More generally, for any operator $\operatorname T$ as in (\ref{e.U}), we can 
  make the operator larger by moving all maximal functions to the 
  right of all square functions.  Thus, it suffices to bound 
  operators of the form 
  \begin{equation*}
  \operatorname T=\operatorname {S}_1\cdots \operatorname{S}_v \operatorname {M}_{v+1}\cdots \operatorname{M}_d.
  \end{equation*}
  
  Recall that the maximal function is bounded as a vector valued map from $L^p(\ell^2)$ into itself for all 
$1<p<\infty$.  This is a well known result of C.~Fefferman and Stein \cite{MR0284802}.  Namely, we have the estimate 
\begin{equation*}
\NOrm \Bigl[\sum _{n} \abs{ \operatorname {M\cdots M} g_n}^2\Bigr] ^{1/2}.p. {}\lesssim{} \NOrm 
\Bigl[\sum _{n} \abs{  g_n}^2\Bigr] ^{1/2}.p. .
\end{equation*}

For a function $f\in L^p(\mathbb R^d)$ and a dyadic  rectangle $R_v\in \mathcal D ^v$, set 
\begin{align*}
f _{R_v}(x_1 ,\dotsc ,x_v,&x_{v+1} ,\dotsc ,x_d) 
\\&{}
\eqdef{} \varphi _{R_v}(x_1 ,\dotsc ,x_v)
\int f(x_1 ,\dotsc ,x_v,x_{v+1},\dotsc ,x_d) \overline{\varphi _{R_v}(x_1 ,\dotsc ,x_v ) }\; d x_1\dots  d x_v .
\end{align*}
It is a consequence of the one variable Littlewood-Paley inequalities that we have 
\begin{equation*}
\NOrm \Bigl[ \sum _{R_v\in \mathcal D^v} \abs{ f _{R_v}}^2\Bigr] ^{1/2}.p. {}\lesssim{}\norm f .p. .
\end{equation*}
 
To conclude, observe that we have 
\begin{equation*}
\operatorname {T}f {}\lesssim{}  \Bigl[ \sum _{R_v\in \mathcal D^v}
\abs{ \operatorname {M}_{v+1}\cdots \operatorname{M} _{d} f _{R_v}}^2\Bigr] ^{1/2} .
\end{equation*}
The $L^p $ norm of the last quantity is clearly bounded by $\norm f.p.$. 
  \end{proof}
    
  The proof of Theorem~\ref{t.2} for a particular range of indices can now 
  be given.  Suppose that we are seeking to bound the paraproduct $\operatorname B$
  from $L^{p_1}\otimes L ^{p_2} $ into $L^r$ where $1<r<\infty$.  Then, 
  by (\ref{e.ipB})
  and the previous lemma, we have 
  \begin{equation*}
  \ip {\operatorname B(f_1,f_2)}, f_3,
  \le{} \prod _{k=1}^3\norm \operatorname T_ k f_k.p_k.
  \end{equation*}
  where by abuse of notation we take $p_3$ to be the conjugate index to $r$.  
  
  The remainder of the theory that we develop is needed to address the case in which 
  the paraproduct does not obey a duality estimate.


\subsection{A Localization Lemma.} 

We will need  an  estimate which refines the $L^2$ estimates for the operators $\operatorname T$
proved in Proposition~\ref{p.Tbounded}. 

We will make a further definition in which these operators, defined as 
as mixture of sums and supremums over rectangles, are restricted to a subset 
of rectangles.  Thus, if $\mathcal O\subset \mathcal R$, and $\operatorname T= \operatorname {SSS }$, 
for instance in the three parameter case, we set 
\begin{equation*}
\operatorname T _{\mathcal O}f=\Bigl[ \sum _{R\in \mathcal O}
\frac{\abs{ \ip f, \varphi_{R},} ^{2} } {\abs{ R}  } \mathbf 1 _{R}\Bigr] ^{1/2}.
\end{equation*}
Of course here we insist that the function $\varphi_R$ have zeros in each coordinate. 
More generally, to define $T _{\mathcal O}$, in the expression (\ref{e.repeatednorms}), we restrict the  
rectangles to be in the collection $\mathcal O$ rather than all possible rectangles.

\begin{lemma}\label{l.local}  Suppose that $\mathcal O\subset \mathcal R$ and that 
there is a constant $\mu>1$ so that for a function $f$,  
\begin{equation} \label{e.localassume}
\operatorname {supp}(f)\cap \mu R=\emptyset, \qquad R \in \mathcal O. 
\end{equation}
Then it is the case that we have 
\begin{equation} \label{e.localconclude}
\norm \operatorname T _{\mathcal O} f .2. {}\lesssim{} \mu ^{-N'}\norm f.2. .
\end{equation}
The exponent $N'$ is a function of only the exponent $N$ in the definition of adapted,
 (\ref{e.adapted}).
\end{lemma}
 
\begin{proof}
This lemma is a corollary to the proof of $L^2$ boundedness of the operators $\operatorname T_j$, 
and to deduce it, we will rely upon a degree of flexibility built into the definition of 
adapted.

If we knew that the functions $\varphi_R$ were supported on, e.g., ~$\tfrac \mu2 R$, then the conclusion 
of the Lemma would be obvious.  Thus, the problem at hand is one of Schwartz tails.

We make a further specification of the definition of adapted, (\ref{e.adapted}), which is
applied to functions on the real line. 
Fix a constant $K$ and an integer $N$ say that $\varphi$ is \emph{ $(K,N)$--adapted } to an interval $I$ iff 
\begin{equation*} 
 \Abs{\operatorname D ^{n} \varphi(x)} {}\le{}K\abs{I}^{-n-\frac12} \Bigl( \tfrac12+\tfrac {\abs{x-c(I)}}{\abs{I}}\Bigr)^{-N}, 
 \qquad n=0,1. 
 \end{equation*}
In addition, the $L^2$ norm of $\varphi$ is at most one.

 Say that  $\varphi$ is 
 \emph{ $(K,N)$--adapted } to a rectangle $R=\otimes _{j=1}^d R _{(j)}\subset \mathbb R^d $ iff 
$\varphi$
is a product 
 \begin{equation*}
 \varphi(x_1,\dotsc,x_d)=\prod _{j=1}^d \varphi_j(x_j)
 \end{equation*}
 with each $\varphi_j$ { $(K,N)$--adapted } to $R _{(j)}$.  This definition naturally extends to 
 collections of rectangles.  
 
 Now, fix $K_0$ and $N_0$ so that Proposition~\ref{p.Tbounded} holds for all functions 
$\{\varphi_R\mid R\in \mathcal R\}$ which are  $(K_0,N_0)$--adapted to $\mathcal R$
(and have zeros in the right coordinates).  Here, we take $K_0\ge1$ so that there will be 
no difficulties with having the functions $\varphi_R$ be of $L^2$ norm one.

For $N'$ as in the conclusion of the Lemma, set $N_1=N_0+N'$.  Consider functions $\{\varphi_R\mid R\in \mathcal R\}$ 
that are $(K_0,N_1)$--adapted to $\mathcal R$. And let $\operatorname T$ be the corresponding operator constructed 
from these functions.

Observe that we can define a new set of functions 
$\{\widetilde \varphi_R\mid R\in \mathcal R\}$ that are $(2K_0,N_0)$ adapted to $\mathcal R$, and satisfy 
\begin{gather*}
\widetilde \varphi_R (x)=\mu^{N'}\varphi_R(x),\qquad x\not\in  \mu R.
\end{gather*}
In those coordinates $1\le{}j\le{}d$ where there is no zero, 
this is accomplished by multiplying by a smooth function that is zero on a large neighborhood of $R$, and 
and identically $\mu ^{N'/d}$ in $\mathbb R-\mu R _{(j)}$.  If the coordinate has a zero, 
observe that 
\begin{equation*}
\ABs{ \int _{\mathbb R- \mu R_{(j) }} \varphi_{R _{(j)}}(x_j)\; dx_j}\le2K_0 \abs{R _{(j) }} ^{1/2}(\mu) ^{-N_1+1}
\end{equation*}
provided $N_1>2$.
And so we can set $\widetilde\varphi_{R _{(j)}} $ in a neighborhood of $R$ to cancel out this integral.

The operator $\widetilde {\operatorname T}$ constructed from the functions $\{\widetilde \varphi_R\mid R\in \mathcal R\}$
will satisfy an $L^2$ bound that is independent of $\mu$.  
For a function $f$ as in (\ref{e.localassume}), we have 
\begin{equation*}
\mu ^{N'}\operatorname T _{\mathcal O}f=\widetilde {\operatorname T} _{\mathcal O} f.
\end{equation*}
And the right hand side admits an $L^2$ bound independent of $\mu$ and $N'$, so the proof is complete.

\end{proof}

We will also need the following corollary to the previous Lemma. 

\begin{corollary}\label{c.local}  Let $\mathcal O$ be a collection of rectangles whose 
shadow has finite measure.  If $f$ is a bounded function, we have the estimate 
\begin{equation}\label{e.localO}
\norm T _{\mathcal O} f.2. {}\lesssim{} \abs{ \operatorname {sh}(\mathcal O)} ^{1/2}\norm f.\infty..
\end{equation}
\end{corollary}

\begin{proof}
Let $f\in L^\infty$ be bounded by one, 
set $U=\operatorname {sh}(\mathcal O)$, and define $f_k$, for $k\ge0$ as in (\ref{e.inductive}). 
We shall see that 
\begin{equation}\label{e.yy}
\sum _{k} \norm T  _{\mathcal O} f_k.2. {}\lesssim{} 
\abs{ \operatorname {sh}(\mathcal O)} ^{1/2}. 
\end{equation}

Indeed, applying  Lemma~\ref{l.local}, we see  that 
\begin{equation*}
\norm \operatorname T _{\mathcal O} f_k.2. {}\lesssim{}2 ^{-N' k} 
\abs{ \operatorname {sh}(\mathcal O)} ^{1/2}
\end{equation*}
where we can assume that $N'>4$ say. 
\end{proof}

\subsection{The Proof of Theorem~\ref{t.2}} We only treat the bilinear case of the theorem, as the higher order 
linearities are easy to accommodate into this proof.  We also restrict our attention to the two parameter
setting.  Straightforward modifications adapt the argument to an arbitrary number of parameters.
The first cases that we consider 
are those in which $\operatorname B$ is to be mapped into a space $L^r$ with $\frac 12\le{}r\le1$. 

Some of the generalities of the proof of the one dimensional case remain in 
force in the current setting, in particular, it will suffice for us to establish 
(\ref{e.1-Lambda-1}).  That is, we shall demonstrate this:  For all $f_j\in L^{p_j} $ of norm 1 and 
set $E_3\subset \mathbb R^2 $ of measure one, 
there is an open  subset $E_3'\subset E_3 $ of measure at least $\frac 12 $, so that for $f_3 $ a smooth
function compactly supported in ${E_3' } $ and with $L^\infty $ norm at most one,  we have 
\begin{equation}  \label{e.2-2do-restricted-weak}
\sum_R 
\abs R ^{-1/2} \prod_{j=1}^3 \abs{\ip f_j , {\varphi_{j,R} } , } {}\lesssim{}1 .
\end{equation}
Moreover, it suffices to take $f_j$ in a dense class of functions, and so we take $f_1$ and
$f_2$ to be smooth and compactly supported.

Observe that as all $f_j $ smooth and compactly supported, the sum above is at most $1 $ if the sum is 
restricted to rectangles that have 
at least one side length either small or large, as in these cases the inner products above decay  
rapidly.\footnote{With the precise definition of small and large depending upon 
the  functions $f_j $.}  Thus, we can assume that the sum is restricted to a finite number of rectangles, and 
we should provide an estimate for the sum that is independent of the exact number or nature of the rectangles.

This last sum is over positive summands. It will be useful to us to organize the sum over appropriate 
subcollections of $\mathcal R$. 
For a collection of dyadic rectangles $\mathcal O$ set 
\begin{equation}  \label{e.sumO}
\operatorname {Sum} (\mathcal O){}\eqdef{}\sum _{R\in \mathcal O}
\abs{ R}\prod _{j=1}^3\frac { \abs{ \ip f_j,\varphi _{j,R},} }{\sqrt {\abs R}}
\end{equation}
where the functions $\varphi_{j,R}$ are associated with the bilinear paraproduct 
that we considering.  

We will be working with different collections of rectangles $\mathcal O$.  The \emph{shadow} 
of $\mathcal O$ is defined to be 
\begin{equation*}
\operatorname {sh}(\mathcal O) {}\eqdef{}\bigcup _{R\in \mathcal O} R. 
\end{equation*}

\subsubsection*{The Definition of $E'_3$} 
Let $\operatorname T_j$, for $j=1,2,3$, be the three operators as in (\ref{e.ipB}).  (Though at this point 
the function $f_3$ is not yet specified.) 
For the sake of symmetry, set $\operatorname T_0$ to be the maximal operator in $2$ parameters. 

Define  $4\nu\eqdef \min(p_1,p_2) $,  and set 
\begin{align}   \label{e.2-zW0} 
\Omega_{j,\ell}&{}\eqdef{} \{ \operatorname {T}_jf_j>\kappa 2^{\ell}\},\qquad \ell\in \mathbb Z, \ j=1,2
	\\   \label{e.zWm}
\Omega_\ell 	& {}\eqdef{}\bigcup _{j=1}^2 \Omega_{j,\ell},
\\	
	\label{e.2-zW^0} 
	\Omega&\eqdef{}\bigcup_{\ell\in \mathbb N} 
	\{ \operatorname{T}_0 \ind {\Omega_{\ell } }>\tfrac1{100}2 ^{-\nu  \ell } \},
\\ \label{e.2-tildezW}
 \widetilde \Omega & {}\eqdef{}\{ \operatorname T_0\ind {\Omega }>\tfrac12\}.
	\end{align}
In these definitions, we  fix a value of  $\kappa{}\simeq{}1$ so that $\abs{ \widetilde \Omega }<\tfrac12 $, and then take 
$E_3'=E_3\cap\widetilde \Omega^c $, so that the measure of this set is at least $\frac12 $. 
This is possible, since we can estimate, using the $L^2$ bound for the
maximal function and the $L^{p_j}$ bounds for the $\operatorname T_j$,
\begin{align*} 
\abs{ \Omega}&\le{}\sum_{\ell\in \mathbb N}  \abs{ \{ \operatorname{T}_0 \ind {\Omega_{\ell } }>2 ^{- \nu \ell } \} }
\\&{}\le{} K_1\sum_{\ell\in \mathbb N}  2 ^{2 \nu \ell}\abs{ \Omega_\ell }
\\&{}\le{}K_2\sum_{\ell\in \mathbb N} \sum_{j=1}^2 \kappa ^{-p_j} 2 ^{(2 \nu-p_j)\ell}.
\end{align*}
And the last sum is less than $\frac18 $ for a fixed $\kappa {}\simeq{}1 $.

\subsubsection*{The Decomposition of $\mathcal R$} \label{s.Rdecomp} 
We decompose the collection of all rectangles.  
 A rectangle $R $ is 
in $\mathcal O_{j,\ell} $ iff $\ell $ is the greatest integer such that 
\begin{equation}  \label{e.2-Out-j,ell}
\abs{ R\cap \Omega _{j,\ell } }=\abs{ R\cap \{ \operatorname {T}_jf_j>\kappa 2^{\ell}\}}\ge{}\tfrac1 {100 } \abs{ R} . 
\end{equation}
Here, we are extending the definition in (\ref{e.2-zW0}) to $j=1,2,3$ and to $\ell\in \mathbb Z$.

Observe that as each $f_j$ is smooth, it is necessarily the case that $\operatorname T_j f_j$
is a bounded function.  Thus, the definition above makes sense, and for each $j$, 
every rectangle $R\in \mathcal R$ will be a member of some $\mathcal O _{ j,\ell } $, 
for $\ell\in \mathbb Z $.

For integers $\vec\ell=(\ell_1,\ell_2,\ell_3)\in \mathbb O {}\eqdef{}(-\mathbb N)\times (-\mathbb N)\times \mathbb Z $, we define 
\begin{equation} \label{e.2Odef}
\mathcal O _{\vec\ell }\eqdef{}\bigcap _{j=1}^3 \mathcal O _{j,\ell_j } .
\end{equation}
This is not a complete decomposition of the collection of rectangles, a point we return to below. 

We appeal to the principal technical estimate, proved below. 
Observe that, by Lemma~\ref{l.main2tech} and using the notation of (\ref{e.sumO}), it is the case that 
\begin{equation}  \label{e.SumOjl}
\operatorname {Sum}(\mathcal O _{\vec\ell} ){}\lesssim{} 
2 ^{\ell_1+\ell_2+\ell_3}\abs{ \operatorname {sh}(\mathcal O _{\vec\ell}) } .
\end{equation}

We therefore need effective estimates for the shadow.  Since each $\operatorname T_j$ operator 
is bounded on all $L^p$ spaces, we of course have the estimate 
\begin{equation}  \label{e.shadowgeneral}
\begin{split}
\abs{ \operatorname {sh}(\mathcal O _{\vec\ell}) }
&{}\lesssim{} \min_j 2 ^{-p_j \ell_j}
\\ & {}\lesssim{} 2 ^{-\theta_1p_1\ell_1-\theta_2p_2\ell_2-\theta_3p_3\ell_3} ,\qquad \vec\ell\in \mathbb O .
\end{split}
\end{equation}
Here,  $\ell_1,\ell_2<0$ while $\ell_3\in \mathbb Z$, 
and  $\theta_j\ge0$ with $\theta_1+\theta_2+\theta_3=1$. 
Recall that $p_1$ and $p_2$ are specified to us in advance, but as $f_3$ is a bounded 
function on a set of finite measure, we are free to take any value of $1<p_3<\infty$ 
that we wish.  In particular, it is effective to take $p_3$ to be relatively close 
to one for $\ell_3>0$ while we take $p_3$ large for $\ell_3\le0$. 

The sums are treated separately based on the sign of the last coordinate of $\vec\ell\in \mathbb O$.  Combining 
(\ref{e.SumOjl}) and (\ref{e.shadowgeneral}), we see that 
 \begin{equation} \label{e.twofinalee}
 \sum _{\vec\ell\in (-\mathbb N)^3 }\operatorname {Sum}(\mathcal O _{\vec\ell} ) {}\lesssim{} 
  \sum _{\vec\ell\in (-\mathbb N)^3 } 2 ^{ \ell_1 (1-p_1 \theta_1 )+{} 
 \ell_2 (1-p_2 \theta_2 )+{}\ell_3 (1-p_3 \theta_3 ) }.
 \end{equation}
 We should choose $0<\theta_1<\frac 1 {p_1}$, and $0<\theta_2<\frac 1 {p_2}$
 so that $\theta_1+\theta_2<1$.  We are then still free to chose $p_3>1$, but 
 close enough to one so that $p_3 \theta_3=p_3(1-\theta_1-\theta_2)<1$.  Thus, 
 this last sum is no more than a constant.

Let us consider the case of $\ell_1,\ell_2\le0$ while $\ell_3>0$.
 The minimum in (\ref{e.shadowgeneral}) occurs for $j=3$, and we have the estimate 
 \begin{equation} \label{e.twofinale}
 \sum _{\vec\ell\in (-\mathbb N)^2\otimes \mathbb N }\operatorname {Sum}
 (\mathcal O _{\vec\ell} ) {}\lesssim{} 2 ^{\ell_1+\ell_2-\ell_3(p_3-1)}.
 \end{equation}
 This clearly sums to a constant for $p_3 $ sufficiently large.

\bigskip 

Some rectangles are not in the classes defined above. 
To treat the remaining cases, set 
\begin{equation} \label{e.Ps}
\begin{split}
\mathcal P _{\vec \ell}  &{}\eqdef{}\mathcal O _{1,\ell_1}\cap \mathcal O _{2,\ell_2}, 
\\
\vec \ell=(\ell_1,\ell_2)\in \mathbb P &{}\eqdef{}\{\mathbb Z-(-\mathbb N)\}\times \{\mathbb Z-(-\mathbb N)\} .
\end{split}\end{equation} 
This decomposition does not take the role of $f_3$ into account, and so 
our next steps are to deduce information about this function.

Suppose that $R\in \mathcal P _{\vec\ell}$. Then either $\ell_1$ or $\ell_2$ must 
be positive.  Suppose that $\ell_1$ is.  
By (\ref{e.2-zW^0}), it is then the case that, by definition,  
$R\subset \Omega$, but moreover	
\begin{equation*}
2 ^{\nu' \ell_1}R\cap E'_3=\emptyset 
\end{equation*}
where $2 \nu'=\nu$.

Now, the function 
$f_3$ satisfies the conditions of (\ref{e.localassume}), with $\mathcal O= 
\mathcal P _{\vec\ell}$ and 
 $\mu=2 ^{\nu'\ell_1 /2}$. 
Then by Lemma~\ref{l.local}, we have 
\begin{equation}\label{e.adaptedN}
\norm \operatorname T_{3,\mathcal P _{\vec\ell} } f_3.2. {}\lesssim{}2 ^{-\nu N'\ell_1/2}
{}\lesssim{}2 ^{-10\ell_1}
\end{equation}
for appropriate choice of $N$ in (\ref{e.adapted}).

Note that we have proved the inequality 
\begin{equation*}
\norm \operatorname T_{3,\mathcal P _{\vec\ell} } f_3.2. {}\lesssim{}\min (2 ^{-10\ell_1},2 ^{-10\ell_2}).
\end{equation*}
We can then apply (\ref{e.main2techconclude2}) to see that 
\begin{equation*}
\operatorname {Sum}(\mathcal P _{\vec\ell}) 
{}\lesssim{}2 ^{\ell_1+\ell_2}\min (2 ^{-10\ell_1},2 ^{-10\ell_2}).
\end{equation*}
This is clearly summable to a constant over the indices $\mathbb P$, as either 
$\ell_1$ or $\ell_2$ must be positive. 

\smallskip 

This proof will permit e.g.~$p_1=1 $, with the additional hypothesis that 
$ f_1\in H^1$, and that all functions $ \varphi _{1,R}$ satisfy (\ref{e.2-zeros}) 
for all coordinates $1\le  j\le d$.  
By duality, this implies the $ \textup{BMO}$ estimate of our Theorem.

\subsection*{The Endpoint Estimates}

The endpoint estimates concern the case when, say, $p_2=\infty$, which 
is a case not handled in the discussion above.   (Note that assuming that $f_2\in L^\infty $, 
we do not need to make additional assumptions about the zeros of the functions
 $ \varphi _{2,R}$.)

We again  prove (\ref{e.2-2do-restricted-weak}).  And the method of proof is quite 
close to the argument above.  Use the same notation as in (\ref{e.2-zW0}), but 
now define 
\begin{align}
	\label{e.zWend}
	\Omega&\eqdef{}\bigcup_{\ell\in \mathbb N} 
	\{ \operatorname{T}_0 \ind {\Omega_{1,\ell } }>\tfrac1{100}2 ^{-\nu  \ell } \},
\\ \label{e.tildezWend}
 \widetilde \Omega & {}\eqdef{}\{ \operatorname T_0\ind {\Omega }>\tfrac12\}.
\end{align}
We take $E_3'=E_3\cap \Omega$, so that again we have $\abs{ E_3'}\ge\frac12$.

We define the sets $\mathcal O _{j,\ell}$ as in (\ref{e.2-Out-j,ell})
but we shall only use this for $j=1,3$.  Set (in contrast to (\ref{e.2Odef})), 
\begin{equation*}
\mathcal O _{\vec \ell }=\mathcal O _{1,\ell_1}\cap \mathcal O _{3,\ell_3}, 
\qquad \vec \ell=(\ell_1,\ell_3). 
\end{equation*}

We then have the estimate below, as a consequence of Corollary~\ref{c.local} and 
(\ref{e.main2techconclude2}), 
\begin{align*}
\operatorname {Sum}(\mathcal O _{\vec\ell} )&{}\le{}
2 ^{-\ell_1-\ell_3} \abs{ \operatorname {sh}(\mathcal O _{\vec\ell}) } ^{1/2} 
\norm \operatorname T _{\mathcal O _{\vec\ell} } f_2.2.
\\& {}\lesssim{}
2 ^{\ell_1+\ell_3}\abs{ \operatorname {sh}(\mathcal O _{\vec\ell}) } .
\end{align*}

We estimate the shadow.  
\begin{equation*}
\abs{\operatorname {sh}(\mathcal O _{\vec\ell}) } {}\lesssim{}\min( 2 ^{-p_1\ell_1},2 ^{-p_3\ell_3})
\end{equation*}
This follows on the one hand from  the assumption that $f_1\in L^{p_1}$.
But recall that we can choose $1<p_3<\infty$ in an arbitrary fashion, as $f_3$ is bounded 
by one and supported on a set of measure at most one.

Pulling these estimate together, we see that 
\begin{equation*}
\operatorname {Sum}(\mathcal O _{\vec\ell} ){}\lesssim{} 2 ^{  \ell_1(1-\theta_1p_1)+\ell_3(1-\theta_3p_3)},
\end{equation*}
where the $\theta_1, \theta_3$ are non-negative and sum to one.  
The index $p_1$ is specified to us, 
but   $p_3$ can be taken arbitrarily.   For $\vec\ell\in (- \mathbb N)\times   \mathbb Z$, 
we should take $p_3$ close to one for $\ell_3\le0$, but $p_3=4$, say, for $\ell_3>0$. 
Doing so we see that 
\begin{equation*}
\sum _{\vec\ell\in (- \mathbb N)\times \mathbb Z}
\operatorname {Sum}(\mathcal O _{\vec\ell} ) {}\lesssim{}1. 
\end{equation*}

\bigskip 

We turn to the case where $\ell_1>0$.  As before, we should now 
gain additional information about the function $f_3$.  But the reasoning of the 
previous section, and in particular (\ref{e.localconclude}) and (\ref{e.adaptedN}),
leads us immediately to 
\begin{equation*}
\norm \operatorname T _{3,\mathcal O _{1,\ell_1}} f_3 .2. {}\lesssim{}2 ^{-10p_1\ell_1}.
\end{equation*}
On the other hand, Corollary~\ref{c.local} implies that 
\begin{equation*}
\norm \operatorname T _{2,\mathcal O _{1,\ell_1}} f_2 .2. {}\lesssim{}
\abs{ \operatorname {sh}(\mathcal O _{1,\ell})} ^{1/2} {}\lesssim{}2 ^{-\frac 12 p_1\ell_1}\,.
\end{equation*}

Appealing to (\ref{e.main2techconclude3}), we see that 
\begin{equation*}
\operatorname {Sum}(\mathcal O _{\ell_1}) {}\lesssim{} 2 ^{\ell_1(1-10-p_1/2)},
\end{equation*}
which is clearly summable over $\ell_1>0$.

\subsubsection*{The Principal Technical Estimate} 

In this section we isolate the principal technical estimate in proof of Theorem~\ref{t.2}.

\begin{lemma}\label{l.main2tech}  Suppose that for three constants $0<\lambda_j<\infty$, 
$j=1,2,3$ and a collection of rectangles $\mathcal O$ we have 
\begin{equation}  \label{e.main2tech}
\abs{ R\cap \{ \operatorname T_j f_j>\lambda_j\} }\le{}\tfrac1{100} \abs{ R},
\qquad R\in \mathcal O, \qquad j=1,2,3.
\end{equation}
Then we have the estimate 
\begin{equation}  \label{e.main2techconclude}
\operatorname {Sum}(\mathcal O) {}\lesssim{} \abs{ \operatorname {sh}(\mathcal O)} 
\prod _{j=1}^3 \lambda_j .
\end{equation}
Suppose that (\ref{e.main2tech}) does not hold for $j=3$.  Then we have the estimate 
\begin{equation}  \label{e.main2techconclude2}
\operatorname {Sum}(\mathcal O) {}\lesssim{}\lambda _1 \lambda_2 
\abs{ \operatorname {sh}(\mathcal O)} ^{1/2} \norm \operatorname T_{3,\mathcal O} f_3 .2. .
\end{equation}
Suppose that (\ref{e.main2tech}) does not hold for $j=2$ and $j=3$. 
Then we have the estimate 
\begin{equation}
\label{e.main2techconclude3}
\operatorname {Sum}(\mathcal O) {}\lesssim{}\lambda _1 
\norm \operatorname T _{2,\mathcal O} f_2.2. \norm \operatorname T_{3,\mathcal O} f_3 .2. .
\end{equation}
\end{lemma}

We will apply this in settings in which we have a good estimate for 
the shadow of $\mathcal O$ in terms of the $\lambda_j$.

\begin{proof}
Set 
\begin{equation*}
W=\operatorname {sh}(\mathcal O)\cap\bigcap _{j=1}^3 \{\operatorname T_j f_j <\lambda_j\} .
\end{equation*}
Then $R\cap W$ has measure at least $\tfrac{97}{100} \abs{ R}$. 
This permits us to restrict the range of integration below to $W$.
\begin{align*} 
\operatorname {Sum}(\mathcal O)& {}\lesssim{} \int _{W}  \sum_{R\in\mathcal O }  
	 \prod _{j=1}^3 \frac{ \abs{\ip f_j , {\varphi _{j,R } }, } } {\sqrt{\abs R}}  \ind R \; dx
	 \\& {}\lesssim{} \int_{W}   \prod _{j=1}^3 \operatorname T _{j } f_j \; dx 
	 \\& {}\lesssim{}  \abs{ \sh {\mathcal O  }}\prod _{j=1}^3 \lambda_j . 
 \end{align*}
 This proves our first conclusion.  The remaining conclusions follow from the same reasoning, 
 with the use of the Cauchy-Schwarz inequality. 
 \end{proof}

%
%
%
%
%


\end{document}